\documentclass{amsart}       
\usepackage{graphicx}
\usepackage{amssymb}
\usepackage{amsmath}
\usepackage{mathrsfs}

\newtheorem{theorem}{Theorem}
\newtheorem*{theorem*}{Theorem}

\newtheorem*{prop*}{Proposition}

\newtheorem*{definition*}{Definition}

\newtheorem{conjecture}{Conjecture}
\newtheorem*{conjecture*}{Conjecture}

\theoremstyle{definition}
\newtheorem{example}{Example}
\newtheorem{question}{Question}
\newtheorem*{question*}{Question}

\newcommand{\Z}{\ensuremath{\mathbb{Z}}}
\newcommand{\Q}{\ensuremath{\mathbb{Q}}}

\newcommand{\R}{\ensuremath{\mathbb{R}}}
\newcommand{\C}{\ensuremath{\mathbb{C}}}

\newcommand{\A}{\ensuremath{\mathbb{A}}}

\newcommand{\cX}{\ensuremath{\mathscr{X}}}

\newcommand{\cZ}{\ensuremath{\mathscr{Z}}}

\newcommand{\cE}{\ensuremath{\mathscr{E}}}

\newcommand{\Spec}{\ensuremath{\mathrm{Spec}\,}}

\newcommand{\red}{\mathrm{red}}

\newcommand{\an}{\mathrm{an}}

\newcommand{\wt}{\mathrm{wt}}
\newcommand{\lct}{\mathrm{lct}}
\newcommand{\Sk}{\mathrm{Sk}}

\newcommand{\llbr}{[\negthinspace[}
\newcommand{\rrbr}{]\negthinspace]}
\newcommand{\llpar}{(\negthinspace(}
\newcommand{\rrpar}{)\negthinspace)}

\newcommand{\sss}{\subsection*{ }\refstepcounter{subsection}{{\bfseries(\thesubsection)}\ }}
\renewcommand{\thetable}{
\theequation}

\begin{document}

\title{Igusa zeta functions and the non-archimedean SYZ fibration}


\author{Johannes Nicaise}


\address{Imperial College,
Department of Mathematics, South Kensington Campus,
London SW72AZ, UK, and KU Leuven, Department of Mathematics, Celestijnenlaan 200B, 3001 Heverlee, Belgium} \email{j.nicaise@imperial.ac.uk}
\thanks{Johannes Nicaise is supported by the ERC Starting Grant MOTZETA (project 306610) of the European Research Council.}

\maketitle

\begin{abstract}
 We explain the proof, obtained in collaboration with Chenyang Xu, of a 1999 conjecture of Veys about poles of maximal order of Igusa zeta functions.
  The proof technique is based on the Minimal Model Program in birational geometry, but the proof was heavily inspired by ideas coming from non-archimedean geometry and mirror symmetry; we will outline these relations at the end of the paper. This text is intended to be a low-tech introduction to these topics; we only assume that the reader has a basic knowledge of algebraic geometry.
\end{abstract}

\section{Introduction}
Let $f$ be a non-constant polynomial in $\Z[x_1,\ldots,x_n]$, for some $n\geq 1$, with $f(\mathbf{0})=0$. For every prime number $p$, Igusa's $p$-adic zeta function $Z_{f,p}(s)$ of $f$ at $\mathbf{0}$ is a meromorphic function on $\C$ that encodes the numbers of solutions of the congruences $f\equiv 0$ modulo powers of $p$ that reduce to the origin $\mathbf{0}$ modulo $p$. The precise definition will be recalled in Section \ref{sec:igusa}. Igusa has proven that $Z_{f,p}(s)$ is a rational function in $p^{-s}$ by performing a local computation of $Z_{f,p}(s)$ on a log-resolution of $f$ over $\Q$.

This proof provides a relation between the poles of $Z_{f,p}(s)$ and the geometry of a log-resolution; see Section \ref{sec:denef}. However, many fundamental questions about these poles are still unanswered, the principal one being Igusa's Monodromy Conjecture, which relates the poles of $Z_{f,p}(s)$ to roots of the Bernstein polynomial of $f$ -- see \eqref{sss:monconj}. Igusa's proof also implies that the order of a pole of $Z_{f,p}(s)$ is at most $n$, and that, if $p$ is sufficiently large, the real part of a pole is at most $-\lct_{\mathbf{0}}(f)$, the opposite of the log canonical threshold of $f$ at $\mathbf{0}$. In 1999, Veys conjectured that the real part of every pole of order $n$ is equal to $-\lct_{\mathbf{0}}(f)$. This implies, in particular, that such poles of maximal order always satisfy the Monodromy Conjecture.

We have proven Veys's conjecture in collaboration with Chenyang Xu in \cite{NicXu2}. The main lines of the proof are explained in Section \ref{sec:veys}. Although the proof makes heavy use of the Minimal Model Program in birational geometry, an important source of inspiration was provided by a combination of non-archimedean geometry and the theory of mirror symmetry, more precisely Kontsevich and Soibelman's non-archimedean interpretation of the SYZ conjecture. This will be explained in Section \ref{sec:SYZ}.

\subsection*{Acknowledgements}
I would like to thank the organizers of the VIASM Annual Meeting 2017 for the invitation to deliver a lecture at the meeting and to write this survey article for the proceedings. The results presented here are a result of joint work with Mircea Musta\c{t}\u{a} and Chenyang Xu, and it is a pleasure to thank them both for the pleasant and interesting collaboration. I am also indebted to Wim Veys for sharing his ideas on the conjecture that constitutes the main subject of this text.

\section{Igusa's $p$-adic zeta functions}\label{sec:igusa}
\sss We fix a non-constant polynomial $f$ in $\Z[x_1,\ldots,x_n]$, for some $n\geq 1$, with $f(\mathbf{0})=0$. For every prime number $p$ and every integer $d\geq 0$, we set
$$N_{f,p}(d)=\sharp \{x\in (p\Z/p^{d+1}\Z)^n\,|\,f(x)\equiv 0\mod p^{d+1} \}.$$
 Thus $N_{f,p}(d)$ is the number of solutions of the congruence $f\equiv 0$ modulo $p^{d+1}$ that reduce to the origin modulo $p$. To study the asymptotic behaviour of these numbers as $d\to \infty$, we introduce the generating series
 $$P_{f,p}(T)=\sum_{d\geq 0}N_{f,p}(d)T^d\quad \in \Z[[T]].$$
 Now it is natural to ask the following question.

\begin{question}[Borevich-Shafarevich, Problem 9 in Ch.1\S5 of \cite{BorSha}]
Is $P_{f,p}(T)$ a rational function? That is, can we write it as a quotient of two polynomials in $\Z[T]$?
\end{question}

\noindent Recall that the rationality of $P_{f,p}(T)$ is equivalent to the existence of a linear recursion pattern in the sequence of coefficients $(N_{f,p}(d))_{d\geq 0}$.

\sss The series $P_{f,p}(T)$ is easy to compute if the zero locus of $f$ in $\A^n_{\Z}$ is smooth over $\Z$ at the origin of $\A^n_{\mathbb{F}_p}$: by Hensel's lemma,
every solution modulo $p^{d+1}$ lifts to precisely $p^{n-1}$ solutions modulo $p^{d+2}$, for all $d\geq 0$, so that $N_{f,p}(d+1)=p^{n-1}N_{f,p}(d)$. This leads to the formula
 $$P_{f,p}(T)=\frac{1}{1-p^{n-1}T}.$$ However, the computation is substantially more difficult if $f$ has a singularity, already in the case of the plane cusp.

\begin{example}\label{exam:cuspgen}
If  $f=(x_1)^2-(x_2)^3$ and $p>3$, the generating series $P_{f,p}(T)$ is given by the formula
$$P_{f,p}(p^{-2}T)=
\frac{1+(p-1)p^{-1}T+(p-1)p^{-4}T^5-p^{-5}T^6}{(1-p^{-1}T)(1-p^{-5}T^{6})}.$$ 
\end{example}

 The reason for rescaling the variable $T$ by a factor $p^{-2}$ will become apparent below (spoiler: the exponents $5$ and $6$ that appear in the denominator are related to the so-called log canonical threshold $5/6$ of the plane cusp singularity).
 We can expect in general that $P_{f,p}(T)$ will reflect some interesting properties of the singularity of $f$ at $\mathbf{0}$.

 \sss Igusa has proven in \cite{Igu} that $P_{f,p}(T)$ is always rational. A key step in the proof is to rewrite the generating series $P_{f,p}(T)$ as a $p$-adic integral
 $$Z_{f,p}(s)=\int_{(p\Z_p)^n}|f|_p^s|dx|$$ where $s$ is a complex variable and $|dx|$ is the normalized Haar measure on $\Z_p^n$. This $p$-adic integral converges if $\Re(s)>0$ and defines a analytic function on the complex right half-plane, which is called Igusa's $p$-adic zeta function of $f$ at the origin $\mathbf{0}$. It is an easy exercise in $p$-adic integration to express $Z_{f,p}(s)$ in terms of $P_{f,p}(T)$: one has  \begin{equation}\label{eq:tranform}
 P_{f,p}(p^{-n-s})=\frac{1 - p^{n}Z_{f,p}(s) }{1-p^{-s}}.\end{equation}
 Thus it is enough to prove that $Z_{f,p}(s)$ is a rational function in $p^{-s}$ (and, in particular, extends to a meromorphic function on $\C$). Igusa proved this by taking a resolution of singularities for $f$ over $\Q$ and using the change of variables formula for $p$-adic integrals to reduce to the case where $f$ is a monomial, which can be solved by a simple computation.

  \sss Igusa's proof not only establishes the rationality of $Z_{f,p}(s)$ but also provides some interesting information about the possible poles of
$Z_{f,p}(s)$ and their expected orders in terms of the geometry of a resolution of singularities for $f$. This is a striking result, because it relates arithmetic properties of $f$ (the poles of the zeta function) with geometric properties of $f$ (the geometry of a resolution of singularities).
  A completely explicit formula for $Z_{f,p}(s)$ in terms of a resolution of $f$ was later given by Denef, for $p\gg 0$. In the next section, we will review the precise formulations of Igusa's theorem and Denef's formula.

\section{Igusa's theorem and Denef's formula}\label{sec:denef}
\sss First, we need to introduce some notation. Let $h:Y\to \A^n_{\Q}$ be a log-resolution for the morphism $f:\A^n_{\Q}\to \A^1_{\Q}$ defined by the polynomial $f$. This means that  $Y$ is a smooth $\Q$-variety, $h$ is a projective morphism of $\Q$-varieties that is an isomorphism over $\A^n_{\Q}\setminus \mathrm{div}(f)$, and $\mathrm{div}(f\circ h)$ is a strict normal crossings divisor on $Y$. Such a morphism $h$ always exists, by Hironaka's embedded resolution of singularities in characteristic zero.

\sss To every log-resolution $h$, we associate the following numerical invariants. We write
$$\mathrm{div}(f\circ h)=\sum_{i\in I}N_iE_i$$
 where $E_i,\,i\in I$ are the prime components of the divisor $\mathrm{div}(f\circ h)$, and the $N_i$ are their multiplicities. Since $h$ is an isomorphism over the complement of $\mathrm{div}(f)$, we can write the relative canonical divisor of $h$ as
 $$K_{Y/X}=\sum_{i\in I}(\nu_i-1)E_i.$$ The number $\nu_i$ is called the log-discrepancy of $X$ at the divisor $E_i$; these are fundamental invariants in birational geometry. Roughly speaking, the multiplicities $N_i$ measure the complexity of $f$ and the log-discrepancies $\nu_i$ measure the complexity of the resolution $h$.
  For every non-empty subset $J$ of $I$, we set
$$E_J=\cap_{j\in J}E_j,\quad E_J^o=E_J\setminus(\cup_{i\notin J}E_i).$$
 The sets $E_J^o$ form a stratification of $\mathrm{div}(f\circ h)$ into locally closed subsets.

 \begin{theorem}[Igusa \cite{Igu}]
 For every prime number $p$, the zeta function $Z_{f,p}(s)$ lies in the ring $$\Q\left[\frac{1}{p^{as+b}-1}\right]_{a,b\in \Z_{>0}}.$$ If $s_0$ is a pole of order $m$ of $Z_{f,p}(s)$, then there exists a subset $J$ of $I$ of cardinality $m$ such that $E_J^o(\Q_p)\cap h^{-1}(p\Z_p^n)\neq \emptyset$ and $\Re(s_0)=-\nu_j/N_j$ for every $j$ in $J$.
 \end{theorem}

 \begin{theorem}[Denef \cite{Den}]
If the prime number $p$ is sufficiently large, then
$$Z_{f,p}(s)=p^{-n}\sum_{\emptyset \neq J\subset I}\sharp(\overline{E}_J^o(\mathbb{F}_p)\cap \overline{h}^{-1}(\mathbf{0}))\prod_{j\in J}\frac{p-1}{p^{N_js+\nu_j}-1}$$
where $\overline{(\cdot)}$ denotes reduction modulo $p$.
 \end{theorem}
\noindent  To be precise, Denef's formula is valid when the resolution $h$ has ``good reduction modulo $p$'' in a certain technical sense; for our purposes, it suffices to know that this condition is always satisfied for $p\gg 0$.

\begin{figure*}[h!]
  \includegraphics[width=1\textwidth]{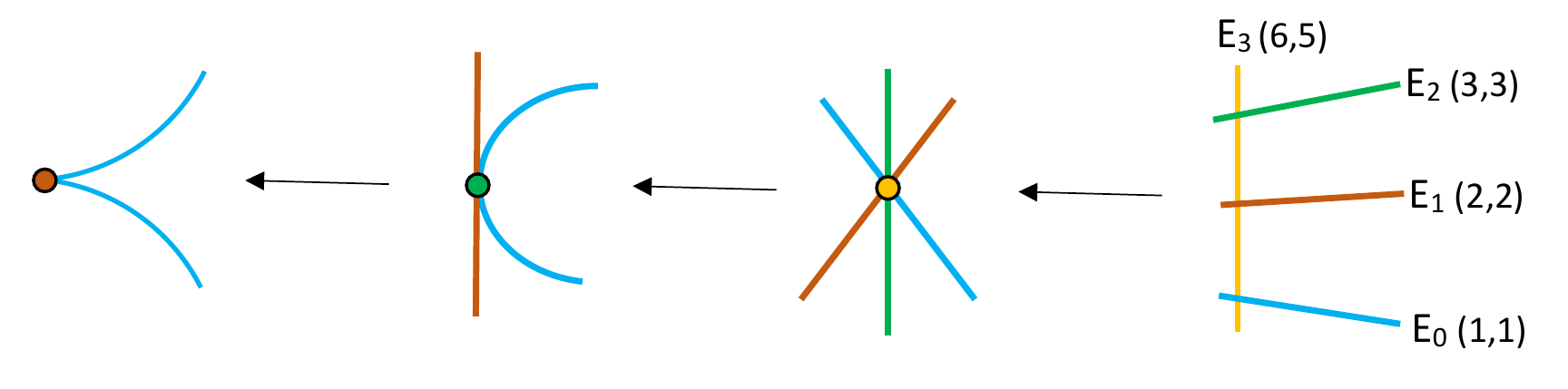}
\caption{Log-resolution of the plane cusp}
\label{fig:cuspres}       
\end{figure*}

\begin{example}\label{exam:cuspzeta}
A log-resolution $h$ for $f=(x_1)^2-(x_2)^3$ can be constructed by means of three consecutive point blow-ups (see Figure \ref{fig:cuspres}). Then $\mathrm{div}(f\circ h)$ has four irreducible components:
 the strict transform $E_0$ of $\mathrm{div}(f)$, with numerical data $N_0=\nu_0=1$, and three exceptional components $E_1,\,E_2,\,E_3$ whose numerical data $(N_i,\nu_i)$ are given by
 $(2,2)$, $(3,3)$ and $(6,5)$, respectively. The log resolution $h$ has good reduction modulo $p$ for all primes $p>3$; in those cases, Denef's formula yields
 \begin{eqnarray*}
 \frac{p^2}{p-1}Z_{f,p}(s)&=& p\frac{p^{-2-2s}}{1-p^{-2-2s}}+p\frac{p^{-3-3s}}{1-p^{-3-3s}}+(p-2)\frac{p^{-5-6s}}{1-p^{-5-6s}}
 \\ & & +(p-1)\frac{p^{-5-6s}}{1-p^{-5-6s}}
 \left(\frac{p^{-1-s}}{1-p^{-1-s}}+\frac{p^{-2-2s}}{1-p^{-2-2s}}+\frac{p^{-3-3s}}{1-p^{-3-3s}}  \right)
 \end{eqnarray*}
 which is compatible with the comparison formula \eqref{eq:tranform} and our expression for $P_{f,p}(p^{-2}T)$ in Example \ref{exam:cuspgen}.
\end{example}

\sss \label{sss:monconj} We can draw some immediate consequences from Igusa and Denef's results. Igusa's theorem implies that the real parts of the poles of $Z_{f,p}(s)$ are all contained
in the finite set
$$\{-\frac{\nu_i}{N_i}\,|\,i\in I,\,E_i(\Q_p)\cap h^{-1}(p\Z_p^n)\neq \emptyset \}.$$
 In practice, most of these candidates will not be real parts of actual poles of $Z_{f,p}(s)$. For one thing, the list of candidates strongly depends on the
 choice of the resolution $h$, whereas $Z_{f,p}(s)$ only depends on $f$ and $p$. But even if $n=2$, when there exists a minimal log-resolution $h$ of $f$, most of
 the candidates will not appear as real parts of poles of the zeta function. A partial explanation of this phenomenon would be given by the so-called Monodromy Conjecture, which puts an additional restriction on the poles of $Z_{f,p}(s)$.

 \begin{conjecture}[Igusa's Monodromy Conjecture]\label{conj:igusa}
If $p$ is sufficiently large and $s_0$ is a pole of $Z_{f,p}(s)$, then $\Re(s_0)$ is a root of the Bernstein polynomial of $f$ at $\mathbf{0}$. In particular, $\exp(2\pi i \Re(s_0))$ is a local monodromy eigenvalue of the complex hypersurface defined by the equation $f=0$.
 \end{conjecture}

\noindent The Monodromy Conjecture quantifies how the singularities of the complex hypersurface defined by $f=0$ influence the poles of the $p$-adic zeta function $Z_{f,p}(s)$, and thus the asymptotic behaviour of the coefficients of the generating series $P_{f,p}(T)$.
 This conjecture has been solved if $n=2$ by Loeser \cite{Loe}, and also for some special classes of singularities, but the general case is wide open. We refer to \cite{Ni:japan} for an overview of results that were known around 2010; a significant breakthrough since that time was the proof for Newton non-degenerate polynomials in three variables -- see \cite{LemPro} and \cite{BorVey}.

\sss Igusa's theorem also implies that the order of a pole $s_0$ of $Z_{f,p}(s)$ is at most
$$\max\{\sharp J\,|\,J\subset I,\,\Re(s_0)=-\nu_j/N_j\mbox{ for all }j\in J,\,E_J(\Q_p)\cap h^{-1}(p\Z_p^n)\neq \emptyset\}.$$
In particular, it is at most $n$, since $E_J$ is empty if $J$ has more than $n$ elements ($n+1$ different prime components of a strict normal crossings divisor on a variety of dimension $n$ can never intersect in a point).

\sss Finally, from Denef's formula, it also follows that the real part of a pole of $Z_{f,p}(s)$ is at most
$$-\min\{\frac{\nu_i}{N_i}\,|\,i\in I,\,E_i\cap h^{-1}(\mathbf{0})\neq \emptyset\}$$
 when $p$ is sufficiently large (if $E_i\cap h^{-1}(\mathbf{0})$ is empty, then $\overline{E}_i\cap \overline{h}^{-1}(\mathbf{0})$ is empty for $p\gg 0$).
 The number
$$\mathrm{lct}_{\mathbf{0}}(f)=\min\{\frac{\nu_i}{N_i}\,|\,i\in I,\,E_i\cap h^{-1}(\mathbf{0})\neq \emptyset\}$$ is an important invariant in birational geometry, called the log-canonical
 threshold of $f$ at $\mathbf{0}$. It is independent of the choice of a log-resolution $h$. It is used to measure the degree of the singularity of $f$ at $\mathbf{0}$, and to divide
 the singularities into different types in the Minimal Model Program.

 \sss The main subject of this survey is the following conjecture.
 \begin{conjecture}[Veys 1999]
 Assume that $p$ is sufficiently large. If $s_0$ is a pole of $Z_{f,p}(s)$ of order $n$, then $\Re(s_0)=-\mathrm{lct}_{\mathbf{0}}(f).$
 \end{conjecture}
\noindent  Thus if  $Z_{f,p}(s)$ has a pole of the largest possible order (namely, $n$), then its real part is also as large as possible. Veys's conjecture was originally stated in \cite{LaeVey} for
  a different type of zeta function (the so-called topological zeta function), but the proof we will present is valid for the topological and motivic zeta functions, as well; see \cite[\S3]{NicXu}. Veys's conjecture implies, in particular, that poles of order $n$ satisfy Igusa's Monodromy Conjecture, because it is known that $-\mathrm{lct}_{\mathbf{0}}(f)$ is the largest root of the Bernstein polynomial of $f$ at $\mathbf{0}$ (see for instance Theorem 10.6 in \cite{Kol}).

\section{Proof of Veys's conjecture}\label{sec:veys}
\sss We will deduce Veys's conjecture from a more general result about the geometry of log-resolutions. Let $k$ be a field of characteristic zero, $X$ a smooth $k$-variety, $f:X\to \A^1_k$ a dominant morphism and $h:Y\to X$ a log-resolution of $f$ as above. We fix a closed point $x$ on $X$ such that $f(x)=0$. The situation we have studied so far corresponds to the case $k=\Q$, $X=\A^n_{\Q}$, $x=\mathbf{0}$. We will continue to use the notations $N_i$, $\nu_i$, $E_J$ etc. The result that we will prove is local on $X$ at the point $x$; shrinking $X$ around $x$, we may assume that $E_J\cap h^{-1}(x)\neq \emptyset$ as soon as $E_J\neq \emptyset$, for every non-empty subset $J$ of $I$. This assumption simplifies some of the notations we will use.

\sss \label{sss:weight} For every $i\in I$, we set $\wt_{f}(E_i)=\nu_i/N_i$ and we call this value the {\em weight} of $f$ at  $E_i$. The log-canonical threshold of $f$ at $x$ is then given by
$$\mathrm{lct}_x(f)=\min\{\wt_f(E_i)\,|\,i\in I\}.$$
We write $Y_0$ for the divisor $$\mathrm{div}(f\circ h)=\sum_{i\in I}N_iE_i,$$
and we denote by $\Delta(Y_0)$ its dual intersection complex. This is a finite $\Delta$-complex whose vertices correspond bijectively to the components $E_i$, $i\in I$, and such that for every subset $J$ of $I$, the number of faces spanned by the vertices $v_j$ with $j\in J$ is equal to the number of connected components of $E_J$.
 We introduce a function
$$\wt_{f}:\Delta(Y_0)\to \R$$ that is completely characterized by the following properties:
\begin{itemize}
\item for every $i\in I$, the value of $\wt_f$ at the vertex of $\Delta(Y_0)$ corresponding to $E_i$ is given by $\wt_f(E_i)=\nu_i/N_i$;
\item the function $\wt_f$ is affine on every face of $\Delta(Y_0)$.
\end{itemize}

\begin{example}\label{exam:weightcusp}
In the case of the minimal log-resolution of the plane cusp (Example \ref{exam:cuspzeta}) we can immediately recover the dual intersection complex $\Delta(Y_0)$ and the
weight function $\wt_f$ from Figure \ref{fig:cuspres}. Since the dimension of $Y_0$ equals one, the dual intersection complex is nothing else than the {\em dual graph} of $Y_0$, with four vertices $v_0,\ldots,v_3$ corresponding to the irreducible components $E_0,\ldots,E_3$ and three edges corresponding to the intersection points of these components. The weight function $\wt_f$ is characterized by its values $\wt_f(v_0)=\wt_f(v_1)=\wt_f(v_2)=1$ and $\wt_f(v_3)=5/6$, and the property that $\wt_f$ is affine on every edge.
\end{example}

\begin{figure*}[h!]
  \includegraphics[width=0.9\textwidth]{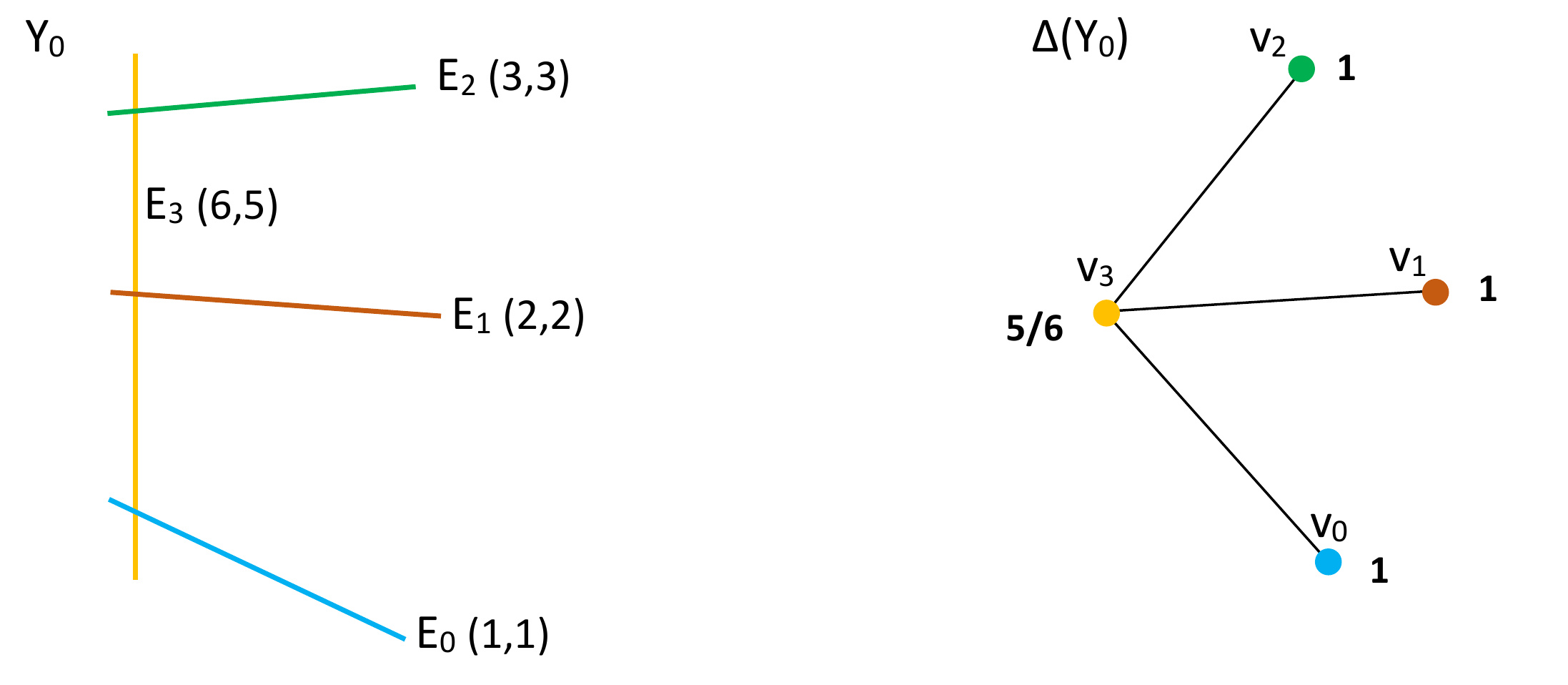}
\caption{The weight function of the minimal log-resolution of the plane cusp}
\label{fig:cuspweight}       
\end{figure*}

\sss \label{sss:observ} The key observation is that $\wt_f$ seems to induce a {\em flow} on $\Delta(Y_0)$ in the direction of decreasing weights, and
 that this flow collapses $\Delta(Y_0)$ onto the locus where $\wt_f$ is minimal, that is, equal to $\lct_x(f)$. Recall that a {\em collapse} if a particular type of
 strong deformation retraction that can be defined in a combinatorial way for $\Delta$-complexes.
 This observation can be made precise in the case $n=2$ (see for instance Example \ref{exam:weightcusp} and Figure \ref{fig:cuspweight}) but it is not clear what the correct notion of a flow
 on a $\Delta$-complex would be in higher dimensions.  We will explain the origin of this observation in Section \ref{sec:SYZ}: it has roots in Kontsevich and Soibelman's non-archimedean SYZ fibration for degenerations of Calabi-Yau varieties in the theory of mirror symmetry.

\sss   In \cite{NicXu}, Chenyang Xu and the author have proven two results that provide compelling evidence for the above observation.
  We denote by $\Delta(Y_0)^{=\lct_x(f)}$ the locus where $\wt_f$ reaches its minimal value $\lct_x(f)$; this is a union of faces of $\Delta(Y_0)$.
 The first result can be understood as follows: if we place ourselves at any point of $\Delta(Y_0)$, then the weight function should somehow select a path along which the weight decreases and which ends in $\Delta(Y_0)^{=\lct_x(f)}$, the locus of minimal weight. If our starting point lies in the interior of a maximal face $\sigma$ of $\Delta(Y_0)$ and the weight function $\wt_f$ is constant on $\sigma$, then $\wt_f$ does not single out any direction at all. So, if our observation is valid, this can only happen when $\sigma$ is already contained in the locus of minimal weight. This property is confirmed by the following theorem.

 \begin{theorem}[Nicaise--Xu, Theorem 2.4 in \cite{NicXu}]\label{thm:maxface}
 If $\wt_{f}$ is constant on a maximal face $\sigma$ of $\Delta(Y_0)$, then its value on $\sigma$ is the minimal value $\mathrm{lct}_x(f)$ of $\wt_{f}$.
\end{theorem}

\noindent The proof of the theorem uses sophisticated techniques from the Minimal Model Program (MMP) in birational geometry; the main idea is to run a well-chosen MMP algorithm on the pair $(Y,Y_0)$ and to study the behaviour of $Y$ near the stratum of $Y_0$ corresponding to the face $\sigma$.

\sss Now we consider the homotopy types of the intermediate level sets of the function $\wt_f$.
  A small technical complication (which will disappear in the SYZ picture in Section \ref{sec:SYZ}) is that the components of $Y_0$ come in two flavours: the exceptional components of $h$ and the components of the strict transform of $\mathrm{div}(f)$. The following theorem refines a fundamental result of de Fernex, Koll\'ar and Xu \cite{dFKoXu}.

 \begin{theorem}[Nicaise--Xu, Theorem 4.10 in \cite{NicXu}]\label{thm:collapse}
 Assume that $\mathrm{div}(f)$ is reduced. We denote by $\Delta(Y_0)_{\mathrm{exc}}$ the subcomplex of $\Delta(Y_0)$ spanned by the vertices that correspond to exceptional components of $h$.
 For every $w\in \R$, we denote by $\Delta(Y_0)^{\leq w}$ the subcomplex of $\Delta(Y_0)$ spanned by the vertices $v$ such that $\wt_{f}(v)\leq w$, and we define $\Delta(Y_0)^{\leq w}_{\mathrm{exc}}$ accordingly.
  \begin{itemize}
 \item Assume that $\lct_x(f)=1$ (that is, $(X,\mathrm{div}(f))$ is log-canonical at $x$). There exists a collapse of $\Delta(Y_0)$ onto $\Delta(Y_0)^{=\lct_x(f)}$ that simultaneously collapses $\Delta(Y_0)^{\leq w}$ onto $\Delta(Y_0)^{=\lct_x(f)}$ for every $w\geq \lct_x(f)$. In particular, the inclusion maps of $\Delta(Y_0)^{=\lct_x(f)}$
     and $\Delta(Y_0)^{\leq w}$ into $\Delta(Y_0)$ are homotopy equivalences.
  \item Assume that $\lct_x(f)\neq 1$. There exists a collapse of $\Delta(Y_0)_\mathrm{exc}$ onto $\Delta(Y_0)^{=\lct_x(f)}$ that simultaneously collapses $\Delta(Y_0)^{\leq w}_\mathrm{exc}$ onto $\Delta(Y_0)^{=\lct_x(f)}$ for every $w\geq \lct_x(f)$. In particular, the inclusion maps of $\Delta(Y_0)^{=\lct_x(f)}$
     and $\Delta(Y_0)^{\leq w}_\mathrm{exc}$ into $\Delta(Y_0)_\mathrm{exc}$ are homotopy equivalences.
  \end{itemize}
    \end{theorem}

Note that in the second case of Theorem \ref{thm:collapse}, we have $\mathrm{lct}_x(f)<1$ so that $\Delta(Y_0)^{=\lct_x(f)}$ is contained in $\Delta(Y_0)_{\mathrm{exc}}$ (components $E_i$ in the strict transform of $\mathrm{div}(f)$ all have $N_i=\nu_i=1$ because we assumed that $\mathrm{div}(f)$ is reduced).

\sss Now let us explain how Theorem \ref{thm:maxface} implies Veys' conjecture. Let $f$ be a non-constant polynomial in $\Z[x_1,\ldots,x_n]$ and let $h:Y\to \A^n_{\Q}$ be a log-resolution for the morphism $f:\A^n_{\Q}\to \A^1_{\Q}$. Suppose that $s_0$ is a pole of order $n$ of $Z_{f,p}(s)$, for some sufficiently large prime number $p$. Denef's formula then implies that there exists a subset $J$ of $I$ of cardinality $n$ such that $E_J\cap h^{-1}(\mathbf{0})\neq \emptyset$ and $\Re(s_0)=-\nu_j/N_j$ for every $j\in J$.
 Then $E_J$ is a finite set of points, and each of these points corresponds to a face of $\Delta(Y_0)$ of dimension $n-1$ on which $\wt_f$ is constant with value $-\Re(s_0)$. Since the dimension of $\Delta(Y_0)$ is at most $n-1$, such a face is always maximal. Now Theorem \ref{thm:maxface} implies that $\Re(s_0)=-\lct_{\mathbf{0}}(f)$, which confirms Veys' conjecture.

\section{The non-archimedean SYZ fibration}\label{sec:SYZ}
\sss To conclude this survey, we explain how the observation in \eqref{sss:observ} originated from the non-archimedean SYZ fibration of Kontsevich and Soibelman.
  The starting point of this story is the theory of {\em mirror symmetry}, which has grown out of mathematical physics (string theory) and has led to fascinating developments in algebraic geometry. For lack of competence, we will only superficially discuss the physical background. According to string theory, our universe can be described as a 10-dimensional real manifold: the 4 dimensions of space-time and an additional 6 dimensions coming from a complex Calabi-Yau  3-fold. Towards the end of the 80s, physicists realized that the Calabi-Yau manifold was not uniquely determined by the physical theory; rather, Calabi-Yau 3-folds seemed to come in {\em mirror pairs} giving rise to equivalent physical theories. One famous property of mirror 3-folds $(X,\check{X})$ is the following symmetry of their Hodge diamonds:
$$h^{p,q}(X)=h^{3-p,q}(\check{X}).$$ Systematic searches in databases of Calabi-Yau 3-folds supplied empirical evidence for the mirror phenomenon.

\sss Mirror symmetry caught the attention of algebraic geometers through a spectacular application to enumerative geometry by Candelas, de la Ossa, Green and Parkes around 1990.
 They managed to compute numbers of rational curves of fixed degrees on the quintic threefold $X$; these are now known as Gromov-Witten invariants.
 Their calculation relied on the insight that the complex and symplectic geometries of $X$ and its mirror partner $\check{X}$ are swapped: the Gromov-Witten invariants
 of $X$ (which are symplectic in nature) can be recovered from certain period integrals on $\check{X}$ (which live in complex geometry and are easier to compute). At that time,
 there was no mathematical framework available to explain these results; a rigorous mathematical confirmation was given only around 1995.

 \sss An important challenge for algebraic geometers was then to give an exact definition of what it means to be a mirror pair and to devise techniques to construct such pairs.
 In recent years, much progress has been made, in particular by Kontsevich--Soibelman and Gross--Siebert. Both of these programs are based on a conjectural
 geometric explanation of mirror symmetry due to Strominger, Yau and Zaslow, known as the SYZ conjecture (1996).  Let $\mathcal{X}^{\ast}$ be a projective family of $n$-dimensional complex Calabi-Yau manifolds
 over a punctured disk $D^{\ast}$ in $\C$, and assume that this family has semistable reduction over the unpunctured disk $D$ and that it is maximally degenerate. The latter condition means that
 the monodromy transformation on the degree $n$ cohomology of the general fiber $\mathcal{X}_t$ of $\mathcal{X}^{\ast}$ has a Jordan block of rank $n+1$. Then $\mathcal{X}_t$ should admit a fibration $\rho:\mathcal{X}_t\to S$ where $S$ is an $n$-dimensional topological manifold and the fibers of $\rho$ are so-called speical Lagrangian tori in $\mathcal{X}_t$.
   As the parameter $t\in D^{\ast}$ tends to $0$, the fibers of $\rho$ are contracted (in the sense of Gromov-Hasudorff convergence of metric spaces) and we are left with the base $S$.
     The mirror partner of $\mathcal{X}_t$ is then constructed by dualizing the torus fibration $\rho$ over the locus where it is smooth and compactifying the result in a suitable way (this involves deforming the dual fibration by so-called quantum corrections, which is the most subtle part of the picture). We refer to the excellent survey paper \cite{Gro} for a more precise statement and additional background on the SYZ conjecture, as well as the Gross--Siebert program.

 \sss The SYZ conjecture seems very difficult to solve in its original differential-geometric set-up.
 A fundamental insight of Kontsevich and Soibelman in \cite{KoSo} is that a close analog of the SYZ fibration can be found in the world of non-archimedean geometry, more precisely in the context of Berkovich spaces. Here, the base $S$ of the fibration arises as a so-called {\em skeleton} and it is intimately related to the geometry of minimal models of the degeneration $\mathcal{X}^{\ast}$ in the sense of the Minimal Model Program. We will now explain the construction of the non-archimedean SYZ fibration, as well as its relation with the results in Section \ref{sec:veys}. Let us emphasize right away that the non-archimedean SYZ fibration is not merely an analog of the conjectural structure in a different context; it can effectively be used to realize the original goal of defining an constructing mirror partners over the complex numbers, by means of non-archimedean GAGA and algebraization techniques. The case of $K3$ surfaces was discussed in \cite{KoSo}, but substantial technical difficulties are still obstructing this approach in higher dimensions.

 \sss We will work over an algebraic model of the complex disk $D$, namely, the spectrum $\Spec R$ of the ring $R=k\llbr t\rrbr$ of formal power series over an algebraically closed  field $k$ of characteristic zero. We write $K=k\llpar t\rrpar$ for the quotient field of $R$; then $\Spec K$ is the algebraic analog of the punctured disk $D^*$.
  The $t$-adic absolute value on $K$, defined by $|a|=\exp(-\mathrm{ord}_t a)$ for all non-zero elements $a$ in $K$, endows $K$ with an analytic structure and allows one to develop a theory of analytic functions and analytic spaces over $K$. The theory that is most suited for our purposes is Berkovich's theory of $K$-analytic spaces, first presented in \cite{Ber}. To every $K$-scheme of finite type $X$, one can attach a $K$-analytic space $X^{\an}$ which is called the {\em analytification} of $X$, in analogy with the situation over the complex numbers. If $X$ is proper over $K$, then the non-archimedean GAGA principle guarantees that $X^{\an}$ completely determines the geometry of $X$; the upshot is that we can apply additional tools to analyze $X^{\an}$ coming from analytic geometry.

  \sss Now let $Z$ be a smooth, projective, geometrically connected $K$-scheme with trivial canonical line bundle. For instance, if $\cX^{\ast}$ is the family of Calabi-Yau manifolds considered above, we get such an object $Z$ over $K=\C\llpar t\rrpar$ by base change to the quotient field of the completed local ring of $D$ at the origin.
  An {\em snc model} for $Z$ is a regular flat projective $R$-scheme $\cZ$ endowed with an isomorphism of $K$-schemes $\cZ_K\to Z$ such that the special fiber $\cZ_k$ is a divisor with strict normal crossings. Then we can attach to $\cZ_k$ a dual intersection complex $\Delta(\cZ_k)$ in the same way as before. A fundamental theorem of Berkovich implies that there exists a canonical topological embedding of $\Delta(\cZ_k)$ into $Z^{\an}$, whose image is called the {\em skeleton} of $\cZ$ and is denoted by $\Sk(\cZ)$. Moreover, this embedding has a canonical retraction $\rho_{\cZ}:Z^{\an}\to \Sk(\cZ)$ that extends to a strong deformation retraction of $Z^{\an}$ onto $\Sk(\cZ)$. We refer to \cite[\S2.5]{Ni:simons} for an explicit example. Thus $Z^{\an}$ is homotopy equivalent to $\Delta(\cZ_k)$, for every snc model $\cZ$ of $Z$. We can use the retraction maps $\rho_{\cZ}$ to organize the skeleta $\Sk(\cZ)$ into a projective system, indexed by the snc models $\cZ$ ordered by the dominance relation. Then one can show that the natural continuous map
  $$Z^{\an}\to \lim_{\stackrel{\longleftarrow}{\cZ}}\Sk(\cZ) $$ induced by the retractions $\rho_{\cZ}$
  is a homeomorphism. This result gives a complete description of the topology of $Z^{\an}$ in terms of the birational geometry of its snc models. It also provides a natural bridge between birational geometry and non-archimedean geometry.

\sss Let $\omega$ be a volume form on $Z$.
 In a joint work with Mircea Musta\c{t}\u{a} \cite{MuNi}, we defined a {\em weight function}
 $$\wt_{\omega}:X^{\an}\to \R\cup \{+\infty\}$$
  which refines a construction of Kontsevich and Soibelman in \cite{KoSo} and is also related to work of Boucksom, Favre and Jonsson in \cite{BoFaJo}.
  Let $\cZ$ be an snc model of $Z$. Then we can view $\omega$ as a rational section of the relative log canonical line bundle $\omega_{\cZ/R}(\cZ_{k,\red})$; we denote the associated Cartier divisor by $\mathrm{div}_{\cZ}(\omega)$.
    The restriction of $\wt_{\omega}$ to the skeleton $\Sk(\cZ)$ of  $\cZ$ can be computed in the following way. It is affine on every face, so that
  it is completely determined by its values at the vertices of $\Sk(\cZ)$. Each vertex $v$ corresponds to a prime component $E$ of the special fiber $\cZ_k$; the value of $\wt_{\omega}$ at $v$ equals $\nu/N$, where $N$ is the multiplicity of $E$ in $\cZ_k$ and $\nu$ is the multiplicity of $E$ in $\mathrm{div}_{\cZ}(\omega)$.
     Note the similarity with the definition of the weight function $\wt_f$ in \eqref{sss:weight}. The weight function $\wt_f$ can also be interpreted as the weight function of a volume form on a suitable $K$-analytic space; see \cite[\S6]{MuNi}.

\sss
The {\em essential skeleton} $\Sk(Z)$ of $Z$ is the set of points in $Z^{\an}$ where $\wt_{\omega}$ reaches its minimal value. This is a non-empty compact subspace of $Z^{\an}$. One can show that the weight function $\wt_{\omega}$ is strictly decreasing under the retraction map $\rho_{\cZ}$, for every snc model $\cZ$ of $Z$, so that $\Sk(Z)$ is always contained in $\Sk(\cZ)$. More precisely, since $\wt_{\omega}$ is affine on every face of $\Sk(\cZ)$, the essential skeleton $\Sk(Z)$ is a union of faces of $\Sk(\cZ)$. The crucial point is that, whereas $\Sk(\cZ)$ obviously depends on $\cZ$, the essential skeleton $\Sk(Z)$ does not depend on the choice of an snc model. It is also independent of $\omega$, because multiplying $\omega$ with a constant in $K^{\ast}$ shifts the weight function by a constant integer.
 In Kontsevich and Soibelman's vision of the non-archimedean SYZ fibration, the base of the fibration is precisely the essential skeleton $\Sk(Z)$. In general, there does not exist any snc model $\cZ$ of $Z$ such that the essential skeleton $\Sk(Z)$ coincides with $\Sk(\cZ)$. However, it was shown by Chenyang Xu and the author in \cite{NicXu2} that we can identify $\Sk(Z)$ with the skeleton of a suitable model of $Z$ if we enlarge the class of models by considering so-called {\em divisorially log terminal (dlt)} models.
  For such (good) dlt models, one can generalize the construction of the skeleton $\Sk(\cZ)$ (which is still homeomorphic to the dual intersection complex of $\cZ_k$) and the retraction map $\rho_{\cZ}:Z^{\an}\to \Sk(\cZ)$.
  If $\cZ$ is a minimal dlt model of $Z$ in the sense of the Minimal Model Program, then $\Sk(Z)=\Sk(\cZ)$.

  \sss This interpretation of $\Sk(Z)$ has two important advantages: it provides us with a map $\rho_{\cZ}:Z^{\an}\to \Sk(Z)$, which plays the role of the SYZ fibration in the non-archimedean context (beware that it depends on the choice of a minimal dlt model, which is not unique in general). Moreover, we can use tools from the Minimal Model Program to study the properties of the essential skeleton $\Sk(Z)$ and the map $\rho_{\cZ}$. According the the SYZ heuristic, the fibers of the non-archimedean SYZ fibration $\rho_{\cZ}$ should be contracted in the limit $t \to 0$, in a suitable sense, so that
  $Z^{\an}$ collapses onto the locus $\Sk(Z)$ where the weight function $\wt_{\omega}$ reaches its minimal value.
      Then it is tempting to speculate that the weight function $\wt_{\omega}$ induces a flow on $Z^{\an}$ towards $\Sk(Z)$, along paths of decreasing weight.
      It is not yet clear what an appropriate definition of flow would be in this context (except if $\mathrm{dim}(Z)=1$; see Example \ref{exam:syz}).
 However, this heuristic has led us to formulate and prove the following theorems, which at their turn were the source of Theorems \ref{thm:maxface} and \ref{thm:collapse} above.

 \begin{theorem}[Nicaise--Xu, Theorem 5.4 in \cite{NicXu}]\label{thm:maxface2} Let $\cZ$ be an snc model of $Z$.
  If $\wt_{\omega}$ is constant on a maximal face $\sigma$ of $\Sk(\cZ)$, then its value on $\sigma$ is the minimal value of $\wt_{\omega}$. In other words,
 $\sigma$ is contained in the essential skeleton $\Sk(Z)$ of $Z$.
\end{theorem}

 \begin{theorem}[Nicaise--Xu, Theorem 5.6 in \cite{NicXu}]\label{thm:collapse2}
 Denote by $w_{\min}$ the minimal value of $\wt_{\omega}$ on $Z^{\an}$, and let $\cZ$ be an snc model of $Z$.
 For every $w\in \R$, we denote by $\Sk(\cZ)^{\leq w}$ the subcomplex of $\cZ$ spanned by the vertices $v$ such that $\wt_{\omega}(v)\leq w$.
  There exists a collapse of $\Sk(\cZ)$ onto $\Sk(Z)$ that simultaneously collapses $\Sk(\cZ)^{\leq w}$ for every $w\geq w_{\min}$. In particular, the inclusion maps of $\Sk(Z)$ and $\Sk(\cZ)^{\leq w}$ into $\Sk(\cZ)$ are homotopy equivalences.
\end{theorem}

\noindent Combining Theorem \ref{thm:collapse2} with Berkovich's strong deformation retraction of $Z^{\an}$ onto $\Sk(\cZ)$, one can then show that
 for every minimal dlt model $\cZ$ of $Z$, the retraction map $\rho_{\cZ}:Z^{\an}\to \Sk(Z)$ can be extended to a strong deformation retraction of $Z^{\an}$ onto the essential skeleton $\Sk(Z)$, which is in line with the SYZ heuristic.

\sss There are some interesting open questions related to the shape of the essential skeleton $\Sk(Z)$. Assume that $Z$ has semi-stable reduction over $R$ (that is, it has an snc model with reduced special fiber) and that it is maximally degenerate (meaning that the monodromy transformation on the $\ell$-adic cohomology of $Z$ in degree $n=\mathrm{dim}(Z)$ has a Jordan block of size $n+1$). Then it was shown in \cite[4.2.4]{NicXu2} that $\Sk(Z)$ is a closed pseudo-manifold of dimension $n$.
If $Z$ is an abelian variety, then $\Sk(Z)$ is homeomorphic to a real torus of dimension $n$, by \cite[4.3.3]{HalNic}.
 Now assume that $Z$ is Calabi-Yau, in the sense that it is geometrically simply connected and the Hodge numbers $h^{0,i}(Z)$ vanish for
 $0<i<n$. Then $\Sk(Z)$ has the rational homology of the $n$-sphere $S^n$ by \cite[4.2.4]{NicXu2}, and its fundamental group has trivial profinite completion \cite[6.1.3]{HalNic}. It is expected that $\Sk(Z)$ is homeomorphic to $S^n$, but
  this appears to be a very hard question; it has been solved by Koll\'ar and Xu in \cite{KolXu} for $n\geq 3$, and for $n=4$ assuming that $Z$ has a minimal dlt model that is snc.
 A crucial difficulty is proving that $\Sk(Z)$ is a topological manifold.

\begin{example}\label{exam:syz}
The easiest example of a non-archimedean SYZ fibration is provided by a Tate elliptic curve over $K$. This is an elliptic curve $E$ over $K$ with Kodaira--N\'eron reduction type $I_n$, for some $n>0$. Every smooth and proper $K$-curve of positive genus has a unique minimal snc model; the condition of having reduction type $I_n$ with $n\geq 2$ is then equivalent with the property that the special fiber of the minimal snc model $\cE$ of $E$ is a reduced cycle of $n$ rational curves (for simplicity we exclude the case $n=1$: then the minimal {\em normal crossings} model has a reduced special fiber which is a rational curve intersecting itself at one point).

In this case, $\cE$ is also a minimal dlt model of $E$. Thus the essential skeleton $\Sk(E)$ equals the skeleton $\Sk(\cE)$ of $\cE$, which we can identify with the dual graph of $\cE_k$; this is a circle subdivided by $n$ vertices (see Figure \ref{fig:SYZ} for an example where $n=3$).

 To reconstruct the entire $K$-analytic space $E^{\an}$, we take the projective limit of the dual graphs $\Sk(\cE')$ where $\cE'\to \cE$ is a composition of blow-ups of points in the special fiber. The effect of such a point blow-up can take two forms: either the center of the blow-up is an intersection point of components in the special fiber, in which case the corresponding edge of the skeleton is subdivided by adding a vertex in its interior; or the center lies on a unique component of the special fiber, in which case a new edge is attached to the vertex of the dual graph corresponding to this unique component. This is illustrated in Figure \ref{fig:SYZ}.

 \begin{figure*}[h!]
  \includegraphics[width=1\textwidth]{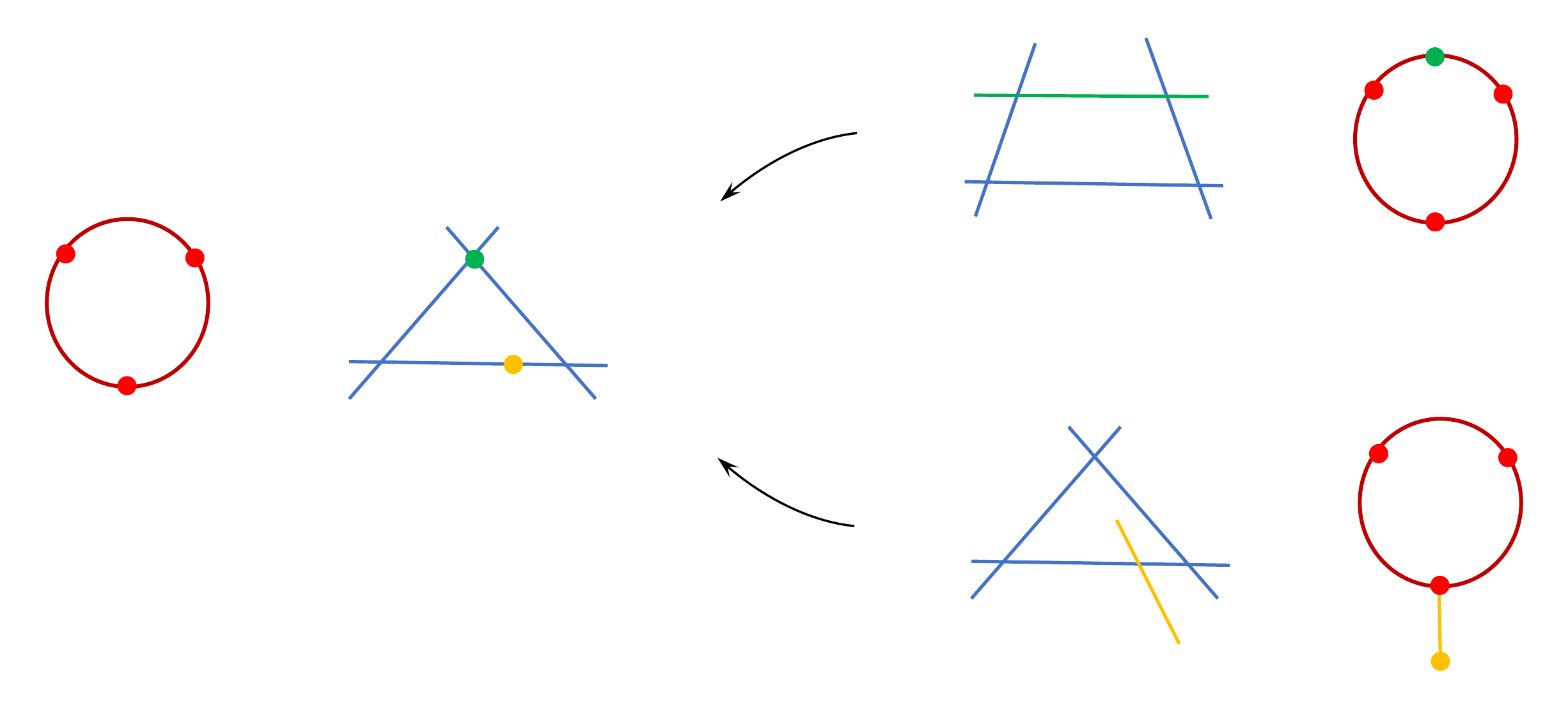}
\caption{Effect of point blow-ups on the dual graphs}
\label{fig:SYZ}       
\end{figure*}

 Repeatedly blowing up points and passing to the projective limit,
 what appears is a circle $\Sk(E)$ with infinitely many trees growing out of it, each tree splitting into infinitely many branches at infinitely many points (Figure \ref{fig:tate}). The non-archimedean SYZ fibration $\rho_{\cE}$ contracts all of these trees onto the circle. For every volume form $\omega$ on $E$, the weight function $\wt_{\omega}$ is constant along the circle $\Sk(E)$ and strictly increasing as we move away from $\Sk(E)$ along one of the branches.

 \begin{figure*}[h!]
\begin{center}
  \includegraphics[width=0.3\textwidth]{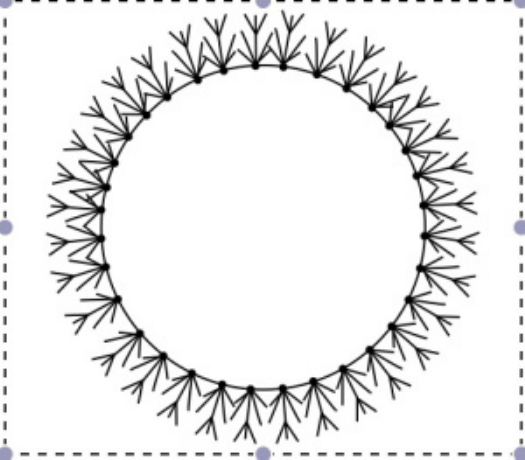}
  \end{center}
\caption{Analytification of a Tate elliptic curve; the essential skeleton is the circle in the middle}
\label{fig:tate}       
\end{figure*}

\end{example}


\begin{thebibliography}{}

\bibitem[Be90]{Ber}
V.~G. Berkovich.
 {\em {Spectral theory and analytic geometry over non-archimedean
  fields}}. Volume~33 of {\em Mathematical Surveys and Monographs}.
 American Mathematical Society, Providence, RI (1990)

\bibitem[BS66]{BorSha}
A.I.~Borevich and I.R.~Shafarevich. {\em Number theory.} Volume 20
of {\em Pure and Applied Mathematics.} Academic Press, New
York-London (1966)

\bibitem[BV16]{BorVey}
B.~Bories and W.~Veys.
 Igusa's $p$-adic local zeta function and the monodromy conjecture for non-degenerate surface singularities.
 {\em Mem. Amer. Math. Soc.} 242  (2016)

\bibitem[BFJ08]{BoFaJo}
S.~Boucksom, C.~Favre and M.~Jonsson. Valuations and
plurisubharmonic singularities.  {\em Publ. Res. Inst. Math. Sci.}
44(2), pages 449--494 (2008)

\bibitem[dFKX17]{dFKoXu}
T.~de Fernex, J.~Koll\'ar and C.~Xu.
  The dual complex of singularities. In: {\em Higher dimensional algebraic geometry, in honour of Professor Yujiro Kawamatas 60th birthday.} Volume 74 of  {\em Adv. Stud. Pure Math.}, pages 103--130. Amer. Math.
Soc., Providence, RI (2017)

\bibitem[De91]{Den}
J.~Denef.  Local zeta functions and Euler
characteristics. {\em Duke Math. J.} 63(3), pages 713--721 (1991)

\bibitem[Gr13]{Gro}
M.~Gross.  Mirror symmetry and the Strominger-Yau-Zaslow
conjecture. In: {\em Current developments in mathematics 2012}, pages 133--191. Int. Press, Somerville, MA (2013)

\bibitem[HN17]{HalNic}
L.H.~Halle and J.~Nicaise. Motivic zeta functions of degenerating Calabi-Yau varieties. Preprint, arXiv:1701.09155.

\bibitem[Ig75]{Igu}
J.~Igusa.\newblock{ Complex powers and asymptotic expansions. II.}
{\em J. Reine Angew. Math.} 278/279, pages 307--321 (1975)

 \bibitem[Ko97]{Kol}
J.~Koll\'ar.
 Singularities of pairs.
In: {\em Algebraic geometry -- Santa Cruz 1995.}
 Volume 62 of {\em Proc. Sympos. Pure Math.}, Part 1, pages 221--287. Amer. Math.
Soc., Providence, RI (1997)

\bibitem[KX16]{KolXu}
J.~Koll\'ar and C.~Xu.
\newblock
{The dual complex of Calabi-Yau pairs.}
\newblock {\em Invent. Math.} 205(3), pages 527--557 (2016)

\bibitem[KS06]{KoSo}
M.~Kontsevich and Y.~Soibelman.  Affine structures and
non-archimedean analytic spaces.
 In: P.~Etingof, V.~Retakh and I.M.~Singer (eds). {\em The unity of
mathematics. In honor of the ninetieth birthday of I. M. Gelfand.}
Volume~244 of {\em Progress in Mathematics}, pages  312--385. Birkh\"{a}user
Boston, Inc., Boston, MA (2006)

\bibitem[LV99]{LaeVey}
A.~Laeremans and W.~Veys.
 On the poles of maximal order of the topological zeta function.
 {\em Bull. London Math. Soc.} 31, pages 441--449 (1999)

\bibitem[LVP11]{LemPro}
A.~Lemahieu and L.~Van Proeyen.
  Monodromy conjecture for nondegenerate surface singularities.
 {\em Trans. Amer. Math. Soc.} 363(9), pages 4801--4829 (2011)

\bibitem[Lo88]{Loe}
F.~Loeser.
 Fonctions d'{I}gusa $p$-adiques et polyn{\^o}mes de {B}ernstein.
 {\em Amer. J. Math.} 110(1), pages 1--21 (1988)

\bibitem[MN15]{MuNi}
M.~Musta\c{t}\u{a} and J.~Nicaise. Weight functions on non-archimedean analytic spaces and the
Kontsevich-Soibelman skeleton.
{\em Alg.~Geom.} 2(3), pages 365--404 (2015)

 \bibitem[Ni10]{Ni:japan}
J.~Nicaise.
 An introduction to $p$-adic and motivic zeta functions and the monodromy conjecture.
 In: G. Bhowmik, K. Matsumoto and H. Tsumura (eds.), {\em Algebraic and analytic aspects of zeta functions and L-functions.}
 Volume 21 of {\em MSJ Memoirs}, pages 115--140. Math.~Soc.~of Japan
 (2010)

\bibitem[Ni16]{Ni:simons}
J.~Nicaise.
 {Berkovich skeleta and birational geometry}. In: M. Baker and S. Payne (eds.), {\em Nonarchimedean and Tropical Geometry}.  Simons Symposia, pages 173--194. Springer, Cham. (2016)



\bibitem[NX16a]{NicXu2}
J.~Nicaise and C.~Xu.
\newblock The essential skeleton of a degeneration of algebraic varieties.
\newblock  {\em Amer.~Math.~J.},  138(6):1645--1667, 2016.

 \bibitem[NX16b]{NicXu}
J.~Nicaise and C.~Xu. Poles of maximal order of motivic zeta functions. {\em Duke Math. J.} 165:2, pages 217--243 (2016)

\end{thebibliography}
\end{document}